\documentclass[12pt]{article}
\usepackage{amsfonts, amsmath, amssymb,amsthm,bbm,comment,fancyhdr,float,makeidx,mathrsfs,verbatim}
\normalsize

\newcommand{\vs}{\vspace*}

\newcommand{\nin}{\noindent}

\newtheorem{mthm}{Theorem}[section]
\newtheorem{mylem}[mthm]{Lemma}
\newtheorem{myprn}[mthm]{Proposition}
\newtheorem{mycor}[mthm]{Corollary}
\newtheorem{mydef}[mthm]{Definition}
\newtheorem{myrem}[mthm]{Remark}
\newtheorem{mycon}[mthm]{Construction}
\newtheorem{myeg} [mthm]{Example}
\newtheorem{myque} [mthm]{Question}
\newenvironment{thm}{\begin{mthm}}{\end{mthm}}
\newenvironment{lem}{\begin{mylem}}{\end{mylem}}

\newenvironment{prop}{\begin{myprn}}{\end{myprn}}

\newenvironment{rem}{\begin{myrem}\rm}{\end{myrem}}

\newenvironment{ex}{\begin{myeg}\rm}{\end{myeg}}
\newenvironment{prof}{\noindent $Proof.$ \rm}{\hfill $\Box$}

\def \nin {\noindent}

\def \Lemma #1 {\vs{3mm}\nin {\bf Lemma #1} \it}
\def \Prop #1 {\vs{3mm}\nin {\bf Proposition #1} \it}
\def \Th #1 {\vs{3mm}\nin {\bf Theorem #1} \it}
\def \Cor #1 {\vs{3mm}\nin {\bf Corollary #1} \it}
\def \Ex #1 {\vs{3mm}\nin {\bf Example #1} \it}

\def \part #1 {\hfil\break\hglue 12pt {\rm (#1)~}}

\def\fs{\footnotesize}

\begin{document}
\title{\Large\bf On monomial ideals whose Lyubeznik  resolution is minimal\thanks{This work is supported by the National Natural Science Foundation of China (No. 11271250).}}
\author{\small\bf {Jin Guo\thanks{Corresponding author. Email: guojinecho@163.com}, \,  Tongsuo Wu\thanks{Email: tswu@sjtu.edu.cn}}\\
{\small \slshape Department of Mathematics,  Shanghai Jiaotong University}\\
{\small \slshape Shanghai,  200240,  China}\\
\small\bf {Houyi Yu\thanks{Email:yhy178@163.com}}\\
{\small \slshape School  of Mathematics and Statistics,  Southwest University}\\
{\small \slshape Chongqing,  400715,  China}\\
}

\date{}
\maketitle

\nin{\bf Abstract.} {\fs For a monomial ideal $I$, let $G(I)$ be its minimal set of monomial generators. If there is a total order on $G(I)$ such that the corresponding Lyubeznik resolution of $I$ is a minimal free resolution of $I$, then $I$ is called a Lyubeznik ideal. In this paper, we characterize the Lyubeznik ideals, and we discover some classes of Lyubeznik ideals. }\\

\nin{\bf 2010 Mathematics Subject Classification:} 13A15; 13D02; 13F20; 05E40\\

\nin{\bf Key Words:} {\fs Lyubeznik ideal, minimal cover, minimal free resolution. }

\section{Introduction}

Let $S=K[x_1, \cdots, x_n]$ be a polynomial ring over a field $K$, and let $I \subseteq S$ be a monomial ideal. One of the most interesting problem is to find an explicit minimal free resolution of $I$ over $S$. Although it is known that every ideal of $S$ has a minimal free resolution (see, e.g.,  Eisenbud \cite{Ei} or Herzog and Hibi \cite{ Herz}), no general description is known, even for monomial ideals.

For some special classes of monomial ideals,  two  ways were developed for finding their minimal free resolutions explicitly. One way is  by investigating various simplicial complexes, e.g., the Taylor resolution, the Scarf complex and the Lyubeznik resolution. For detailed constructions of the three mentioned simplicial complexes, one can refer to the recent work \cite{Mer} by J. Mermin.  Generally speaking, the Taylor resolution is far from being minimal, while the Scarf complex is minimal but is not exact in general. Although the Lyubeznik resolution is also not minimal in general, it is much closer to the minimal free resolution of an ideal compared with the Taylor resolution. By considering rooting maps, Novik  in \cite{No} proved that the Lyubeznik resolution is a minimal free resolution for the matroid ideal of a finite projective space.
Another way is provided by Eliahou and Kervaire \cite{EK} in 1990, in which they constructed  minimal free resolutions of Borel ideals. Some classes of other monomial ideals were given by using this method, such as the linear quotients ideals with decomposition function by Herzog and Takayama \cite{HT}, the edge ideals with some combinatorial conditions by Horwitz \cite{Ho}, and so on.

The main purpose of the present work is to describe a monomial ideal which can be resolved to obtain a minimal Lyubeznik resolution, under a properly chosen total order on the minimal monomial generating set of the ideal. Apparently, which total order should be chosen must be the central topic of the problem.

In section 2, we introduce some new definitions which will be used in the next sections. In section 3, we describe the Lyubeznik ideals by taking advantage of the concept of E-minimal cover (Theorem \ref{Ly ideal}). Then in section 4, we give some important properties, especially, the description about M-minimal complete cover, to help judge a Lyubeznik ideal conveniently (Proposition \ref{m min cover}). In section 5, we use Theorem \ref{Ly ideal} and Proposition \ref{m min cover} to discover several classes of Lyubeznik ideals. In section 6, we show how to examine some special classes of elements  in $G(I)$, i.e., out points, inner points and boundary points of a cover, defined in the next section.

\section{Some definitions}

For a monomial ideal $I$, let $G(I)=\{u_1, u_2, \cdots, u_s\}$ be its minimal set of monomial generators.
For a subset $A$ of $G(I)$, the multidegree of $A$, denoted by $m(A)$,   is the least common multiple of the elements in $A$.
We call a subset $C$ of $G(I)$ a $\it cover$ of a monomial $u \in C$, if $u \mid m(C\setminus \{u\})$, or alternatively we say $C$ covers $u$, denoted by $u \Box C$. The {\it complete cover} induced by a cover $C$, denoted by $\overline{C}$, is the subset of $G(I)$ that contains all of the elements of $G(I)$ which divide $m(C)$. A cover $C$ (of $u$) is called an {\it M-minimal cover} of $G(I)$, if there exists no cover $V$ (of some $v$) whose multidegree $m(V)$ is a proper factor of $m(C)$. A cover $C$ of a monomial $u$ is called an {\it E-minimal cover} of $u$ if no proper subset of $C$ can  cover $u$.


For example, in the monomial ideal $$I=(x^{4}, y^{4}, x^{3}y, xy^{3}, x^{2}y^{2}),$$
the subset $\{x^{4}, y^{4}, x^{2}y^{2}\}$ is an E-minimal cover of $x^{2}y^{2}$, but it is not an M-minimal cover of $G(I)$. Actually, it is easy to see that the multidegree of the cover $\{x^{3}y, xy^{3}, x^{2}y^{2}\}$ is $x^3y^3$, and it properly divides $x^4y^4$, which is the multidegree of $\{x^{4}, y^{4}, x^{2}y^{2}\}$.


Let $C$ be a cover. We call a subset $D$ of $C$ an {\it out set} of $C$, if $m(D)=m(C)$ and the multidegree of any proper subset of $D$ is not equal to $m(D)$. Of course, an out set of a cover may be not unique. A point is called an {\it out point} of $C$ if it is in every out set of $C$ , and we use $\mathcal{O}(C)$ to denote the out points of $C$. We call a point not in any out set of $C$ an {\it inner point} of $C$, and we use $\mathcal{I}(C)$ to denote the inner points of $C$. The other points in $C$ are called {\it boundary points} of $C$, and the set of all boundary points in $C$ is denoted by $\mathcal{B}(C)$. An element $u\in C\setminus \mathcal{I}(C)$ is called {\it exchangeable } in $C$ if for every out set $D$ of $C$ which contains $u$, we have $m((D \setminus \{u\}) \cup \{v\})= m(D)$ for every $v \in C \setminus D$. The set of all exchangeable points in $C$ is denoted by $\mathcal{E}(C)$. It is easy to see from the definition that $\mathcal{E}(C) \subseteq \mathcal{B}(C)$. For example, in the monomial ideal $I=(xy, yz, xz)$, $G(I)=\{xy, yz, xz\}$. Of course, $C=G(I)$ is a complete cover. It is easy to see that $\mathcal{E}(C)= \mathcal{B}(C)=C$.

Note that the above definitions are independent of a total order on $G(I)$.

Let $\prec$ be a total order on $G(I)$, and let $A$ be a subset of $G(I)$. Let $min(A)$ be the least element of $A$ under the total order $\prec$. Let $B$ be another subset of $G(I)$. If $min(A) \prec min(B)$, then we write $A \prec B$. If $A$ has only one element $u$ and $u\prec min(B)$, then we denote $u \prec B$. A set $D$ is  said to be {\it broken} under the total order $\prec$, if there exists an element $u\in G(I)$, such that $u \mid m(D)$ and $u \prec D$. A subset $E$ of $G(I)$ is called {\it preserved}, if no subset of $E$ is broken.

Let $\bigtriangleup_I$ be the full simplex on $G(I)$. For a given total order $\prec$ on $G(I)$, let $L_{(I, \prec)}$ be the following simplicial subcomplex of $\bigtriangleup_I$:
$$L_{(I, \prec)}=\{F\in \bigtriangleup_I \mid min\{u\in G(I) \mid u|m(G)\} \in G \ \ {\rm for \ all}\ G\subseteq F\}.$$
We call $L_{(I, \prec)}$ a {\it Lyubeznik simplicial complex} of $I$.
The following associated algebraic chain complex is proved to be a free resolution of $I$, and is called the {\it Lyubeznik resolution} of $I$ under the total order $\prec$.
$$\mathcal{L}: \cdots \stackrel{\varphi_{n}}{\longrightarrow} L_n \stackrel{\varphi_{n-1}}{\longrightarrow} L_{n-1} \cdots  \stackrel{\varphi_{1}}{\longrightarrow} L_1 \stackrel{\varphi_{0}}{\longrightarrow} I \longrightarrow 0 $$
In this resolution, $L_i=\{F \in L_{(I, \prec)} \mid |F|=i\}$. For a given $F=\{u_{j_1}, u_{j_2}, \cdots, u_{j_i}\} \in L_i$, let $G_k=F\setminus \{u_{j_k}\} \in L_{i-1}, 1\leq k \leq i$. Recall that $\varphi_{i-1}(F)=\sum\limits_{k=1}^i \varepsilon_F^{G_k} \frac{m(F)}{m(G_k)} G_k$, where the sign $\varepsilon_F^{G_k}$ equals to $1$ (respectively, $-1$ ) for odd $k$ (for even $k$, respectively).

Generally speaking, the Lyubeznik resolution $\mathcal{L}$ is not minimal. Let $\bf m$ be the homogeneous maximal ideal of $S$, i.e., $\bf{m}$ $=(x_1, x_2, \cdots, x_n)$. $\mathcal{L}$ is minimal if and only if $\varphi_{i-1}(L_i) \subseteq \bf{m}$$L_{i-1}$ for all $i$. By the construction of $\varphi$, $\mathcal{L}$ is minimal if and only if $m(F)\neq m(G_k)$ for all $F$ and all $k$.
If there is a total order on $G(I)$ such that the corresponding Lyubeznik resolution of $I$ is minimal, then $I$ is called a {\em Lyubeznik ideal}.

Clearly, a Lyubeznik resolution relies heavily on the given total order $\prec$ on $G(I)$. In fact, even if $I$ is a Lyubeznik ideal, the Lyubeznik resolution of $I$ determined by a given total order $\prec$ on $G(I)$ may be not minimal. For example, consider the monomial ideal $I=(x^3, x^2y, y^3, y^2z, z^3)$ of $S=K[x,y,z]$. On $G(I)$ define a total order $\prec$ by $$x^3 \prec x^2y \prec y^3 \prec y^2z \prec z^3.$$ By definition, both $\{x^3, x^2y, y^2z\}$ and $\{x^3, y^2z\}$ are faces of $L_{(I, \prec)}$. So, $\varphi(L_3) \not\subseteq (x, y, z)L_2$, since $m(\{x^3, x^2y, y^2z\})=m(\{x^3, y^2z\})=x^3y^2z$. Hence the Lyubeznik resolution of $I$ by the total order $\prec$ is not minimal. On the other hand, the Lyubeznik resolution of $I$ determined by another total order $\vdash$ is a minimal free resolution of $I$, where $$x^2y \vdash y^2z \vdash x^3 \vdash y^3 \vdash z^3.$$ So, $I$ is a Lyubeznik ideal.


\section{Lyubeznik ideals}

Let $I$ be a monomial ideal. If there is no cover in $G(I)$ (e.g., $G(I)=\{x^2,\,yz,\,y^2\}$), then clearly $I$ is a Lyubeznik ideal. In fact, the Taylor resolution, the Scarf complex and the Lyubeznik resolution of $I$ under any total order on $G(I)$ are identical.

For a given monomial ideal $I$ with covers in $G(I)$, we give the following criterion on judging wether $I$ is a Lyubeznik ideal.

\begin{thm}\label{Ly ideal}
Let $I$ be a monomial ideal in the polynomial ring $K[x_1, x_2, \cdots, x_n]$ with $G(I)=\{u_1, u_2, \cdots, u_s\}$. Then the following statements are equivalent:

$(1)$  $I$ is a Lyubeznik ideal.

$(2)$ There exists a total order $\prec$ on $G(I)$ such that for every element $u$ of $G(I)$ and every E-minimal cover $C$ of $u$, $C$ is not preserved.

$(3)$ There exists a total order $\prec$ on $G(I)$ such that for every element $u$ of $G(I)$ and every E-minimal cover $C$ of $u$, there exist $D\subseteq C$ and $v\notin D$, such that $D\cup \{v\}$ is an E-minimal cover of $v$, satisfying $min(\overline{D} \setminus D) \prec min(D)$.
\end{thm}

\begin{prof}
(1) $\Longleftrightarrow$ (2):If $I$ is not a Lyubeznik ideal, then for every total order $\prec$ on $G(I)$, and the corresponding Lyubeznik simplicial complex $L_{(I, \prec)}$ with differential map $\varphi$, there exists $i< s$, such that $\varphi(L_{i+1}) \not\subseteq (x_1, x_2, \cdots, x_n)L_i$. So, there exist $D \subseteq G(I)$ and $u \notin D$, such that $D, D\cup \{u\} \in L_{(I, \prec)}$ and $m(D)=m(D\cup \{u\})$. By the construction of $L_{(I, \prec)}$, $D\cup \{u\}$ is a cover of $u$ which is preserved. It is easy to see that there exists $E \subseteq D$, such that $E\cup \{u\}$ is an E-minimal cover of $u$, and $E\cup \{u\}$ is preserved. Conversely, if for every order $\prec$ on $G(I)$, there exists an element $u$  and an E-minimal cover $C$ of $u$ such that $C$ is preserved, then we know that $C\setminus \{u\}$ is preserved as a subset of $C$. By the construction of $L_{(I, \prec)}$, both $C$ and $C\setminus \{u\}$ are in $L_{(I, \prec)}$. So, the Lyubeznik resolution of $I$ based on the order $\prec$ is not minimal since the multidegrees of $C$ and $C\setminus \{u\}$ are identical. Hence $I$ is not a Lyubeznik ideal.
\\(2) $\Longrightarrow$ (3): The condition $C$ is not preserved implies that there exists a subset $D$ of $C$ such that some element $\{v\} \notin D$ is covered by $D\cup \{v\}$, and $v$ is less than every element of $D$ in the order $\prec$. It is easy to choose a subset of $D$, say $E$, such that $E\cup \{v\}$ is an E-minimal cover of $v$. Clearly, $v$ is less than every element of $E$ in the order $\prec$. So, $min(\overline{E} \setminus E) \prec  min(E)$.
\\(3) $\Longrightarrow$ (2): Note that $D$ is broken under condition (3), so the result of (2) is clear.
\end{prof}



 \begin{rem}\label{Ly}$(1)$ Theorem \ref{Ly ideal} actually provides a method for judging whether an ideal is a Lyubeznik ideal. First, list all  the elements of $G(I)$, then compute all the E-minimal covers, and finally judge whether a total order exists which fulfills all the relations determined  by all covers. Actually, we only need to check whether such a partial order exists, since every finite partial order can be refined to be a total order, see Lemma \ref{refine}.

$(2)$ The method for checking if a monomial ideal is a Lyubeznik ideal is very useful for monomial ideals with a small number of monomial generators. But when the generators' number is huge, the verifications would be rather complicated. So in the next sections, we will study further properties of Lyubeznik ideal.
\end{rem}

\vs{3mm}In the following example, we show how to use Theorem \ref{Ly ideal} to judge if a monomial ideal is Lyubeznik.
\begin{ex}\label{L ideal}
Consider the monomial ideal $I=(x^3, x^2y, y^3, y^2z, z^3)$ of $S=K[x,y,z]$. In order to judge whether $I$ is a Lyubeznik ideal, we find all the elements with their E-minimal covers, as the following shows:
$$x^2y \Box \{x^3, x^2y, y^2z\},\,\, x^2y \Box \{x^3, x^2y, y^3\},\,\, y^2z \Box \{y^3, y^2z, z^3\}. $$
By Theorem \ref{Ly ideal}, we need to check if the sets $\{x^3, x^2y, y^2z\}$, $\{x^3, x^2y, y^3\}$, $\{y^3, y^2z, z^3\}$ are preserved. For the first set,  only one subset $\{x^3, y^2z\}$ could be broken, and the only element in $G(I)$ could break it is $x^2y$. By Theorem \ref{Ly ideal}, if we have a total order $\prec$ on $G(I)$ such that $I$ is a Lyubeznik ideal, then we must have $x^2y \prec \{x^3, y^2z\}$ in $\prec$. By the same way, we have $x^2y \prec \{x^3, y^3\}$ and $y^2z \prec \{y^3, z^3\}$. Finally, there exists a total order on $G(I)$, say,
$$x^2y \prec y^2z \prec x^3 \prec y^3 \prec z^3, $$
satisfying the three conditions. So, the Lyubeznik resolution under the order $\prec$ is minimal. Hence $I$ is a Lyubeznik ideal.


\end{ex}

\section{Some properties of Lyubeznik ideals}

Let $I$ be a Lyubeznik ideal under the total order $\prec$. Assume further that $u_1 \prec u_2 \prec \cdots \prec u_s$ on $G(I)$. For every element $u_i \in \mathcal{O}(G(I))$, if we change the order, such that $u_i$ to be the largest one and keep the relations of the other elements, then we get a new order $\vdash$ by the following:
$$u_1 \vdash \cdots \vdash u_{i-1} \vdash u_{i+1} \vdash \cdots \vdash u_s \vdash u_i.$$
It is easy to see that the Lyubeznik resolution of $I$ under the total order $\vdash$ is still minimal. In fact, every set broken in $\prec$ must be broken in $\vdash$. By Theorem \ref{Ly ideal}, the more for sets being broken, the closer for $I$ being Lyubeznik. So we have the following conclusion:


\begin{prop}\label{change place}
If I is a Lyubeznik ideal, then there exists a total order defined by the relations $u_{i_1} \prec \cdots \prec u_{i_a} \prec u_{j_1} \prec \cdots \prec u_{j_b},$ where $\mathcal O(G(I))= \{u_{j_1} \cdots u_{j_b}\}$.
\end{prop}

From now on, we  only need to check whether there is this kind of order which fulfills the condition in Theorem \ref{Ly ideal}. In other words, it is not necessary to consider the out points of $G(I)$ when  judging if $I$ is Lyubeznik. We can simply let all the out points larger than the other points of $G(I)$.

When applying Theorem \ref{Ly ideal}, the most complicated part is to check whether a set is preserved. But it is relatively easy for an M-minimal complete cover, as the following result shows.

\begin{prop}\label{m min cover}
Let I be a Lyubeznik ideal under a total order $\prec$. If $C$ is an M-minimal complete cover of $G(I)$, then either $\mathcal{I}(C)\neq \emptyset$ or $\mathcal{E}(C)\neq \emptyset$. Furthermore, exactly one of the following cases occurs:
\\(1) $\mathcal{I}(C)\neq \emptyset$. In this case, we have $\mathcal{I}(C) \prec C\setminus \mathcal{I}(C)$.
\\(2) $\mathcal{I}(C)= \emptyset$ and $\mathcal{E}(C)\neq \emptyset$. In this case, we have $\mathcal{E}(C) \prec C\setminus \mathcal{E}(C)$.
\end{prop}

\begin{prof}
In the case $\mathcal{I}(C)\neq \emptyset$, if the least element $u$ of $C$ is not an inner point, then there exists an out set of $C$, say $D$, such that $u\in D$. We claim that $D \cup \{v\}$ is preserved for every $v \in \mathcal{I}(C)$. In fact, for a subset $E$ of $D \cup \{v\}$, if $u \in E$, clearly $E$ is not broken; if $u \notin E$, then $m(E)$ is less than $m(D)$ since $v \in \mathcal{I}(C)$. By the definition of M-minimal, $E$ is not broken. Hence $D \cup \{v\}$ is a cover of $v$, and $D \cup \{v\}$ is preserved. By Theorem \ref{Ly ideal}, it is a contradiction.
In the case $\mathcal{I}(C)= \emptyset$ and $\mathcal{E}(C) \neq \emptyset$, if the least element $u$ of $C$ is not an exchangeable point, then there exists an out set of $C$, say $F$, such that $u\in F$. We can choose an element $v$ which can not be exchanged with $u$, which implies $m(F\cup \{v\} \setminus \{u\}) \prec m(F)$. In a similar way, we can also get a contradiction. It is easy to see that if an M-minimal cover contains neither inner point nor exchangeable point, then $I$ is not a Lyubeznik ideal.
\end{prof}

\begin{rem}\label{algorithm} Now let us set up a criterion on how to judge whether a monomial ideal $I$ is a Lyubeznik ideal more easier. At the mean time, we also give an algorithm on how to find a total order $\prec$ for a Lyubeznik ideal $I$.
\\(1) Find all the out points of $G(I)$, and let them larger than other  points.
\\(2) List all the other points and their covers, and then compute the multidegrees of the corresponding covers.
\\(3) Pick up all the M-minimal complete covers, collect all the relations by Proposition \ref{m min cover}. If the relations do not satisfy the transitivity law or the anti-symmetry law, then we conclude that $I$ is not a Lyubeznik ideal. Otherwise, go to the next step.
\\(4) Pick up all the other E-minimal covers, use Theorem \ref{Ly ideal} to collect all the relations, and judge if there exists a total order which is coherent with all the relations.
\end{rem}
\begin{rem}\label{M-E} When judging whether an ideal is Lyubeznik, one may get a wrong answer if he only consider the M-minimal complete covers. For example, the ideal
$$I=(x_1^3x_3^3,\, x_2^3x_4^3x_5^3, \,x_1^4, \,x_3^4, \, x_3^2x_4^2, \,x_3^2x_5^2,\, x_3x_4x_5)$$
has two M-minimal complete covers, i.e.,  $\{x_3^2x_4^2, x_3^2x_5^2, x_3x_4x_5\}$ and $\{x_1^3x_3^3, x_1^4, x_3^4\}$. Each of them has  a unique inner point, $x_3x_4x_5$ and $x_1^3x_3^3$, respectively. So, we have relations $x_3x_4x_5 < \{x_3^2x_4^2, x_3^2x_5^2 \}$ and $x_1^3x_3^3 < \{ x_1^4, x_3^4\}$. If we don't consider the other E-minimal covers, we can get a total order $\prec$ defined by $$x_1^3x_3^3 \prec x_3x_4x_5 \prec x_3^2x_4^2 \prec x_3^2x_5^2 \prec x_2^3x_4^3x_5^3 \prec x_1^4 \prec x_3^4,$$ under which the Lyubeznik resolution of $I$ is clearly not minimal. In fact, because $\{x_1^3x_3^3, x_2^3x_4^3x_5^3, x_3x_4x_5\}$ is an E-minimal cover over $x_3x_4x_5$, but $\{x_1^3x_3^3, x_2^3x_4^3x_5^3, x_3x_4x_5\}$ is preserved, hence the Lyubeznik resolution of $I$ under $\prec$ is not minimal.
Note that the ideal is really a Lyubeznik ideal, if we set the total order $\vdash$ to be $$x_3x_4x_5 \vdash x_1^3x_3^3 \vdash x_3^2x_4^2 \vdash x_3^2x_5^2 \vdash x_2^3x_4^3x_5^3 \vdash x_1^4 \vdash x_3^4, $$ which can be obtained after considering further the E-minimal cover $\{x_1^3x_3^3, x_2^3x_4^3x_5^3, x_3x_4x_5\}$ over $x_3x_4x_5$.
\end{rem}

\vs{3mm}Now let us reconsider Example \ref{L ideal}. First, observe $\mathcal{O}G(I)=\{x^3 , y^3 , z^3\}$ and let them larger than all the other points in the forthcoming total order. Then we find a M-minimal cover, $x^2y \Box \{x^3, x^2y, y^2z\}$ with multidegree $x^3y^2z$. (Actually we have three M-minimal covers in $G(I)$, $x^2y \Box \{x^3, x^2y, y^2z\}$ with multidegree $x^3y^2z$, $x^2y \Box \{x^3, x^2y, y^3\}$ with multidegree $x^3y^3$, $y^2z \Box \{y^3, y^2z, z^3\}$ with multidegree $y^3z^3$ respectively, but the last two covers are covered by the out points of $G(I)$, so by proposition \ref{change place}, it is not necessary to consider them here.) So the possible order must fulfill $x^2y < y^2z$. So we can get a good total order $\prec$ defined by $$x^2y \prec y^2z \prec x^3 \prec y^3 \prec z^3.$$
So, $I$ is a Lyubeznik ideal under the order $\prec$, or under the order $\vdash$ defined by $$x^2y \vdash y^2z \vdash y^3 \vdash z^3 \vdash x^3 $$ by Proposition \ref{change place}.

\vs{3mm}Proposition \ref{m min cover} is powerful when searching for counterexamples of Lyubeznik ideals. We show it in the following example.

\begin{ex}\label{not Ly ideal}
The monomial ideal $I=(x^3, x^2y, y^3, yz^2, z^3)$ is  not a Lyubeznik ideal. In fact, it is easy to see that there are four M-minimal complete covers. Consider two of them, i.e.,  $x^2y \Box \{x^3, x^2y, yz^2\}$ and $yz^2 \Box \{x^2y , yz^2, z^3\}$. By Proposition 4.2, if $I$ is Lyubeznik under a total order $\prec$ on $G(I)$, it  must fulfill both $x^2y \prec \{x^3, yz^2\}$ and $yz^2 \prec \{x^2y , z^3\}$. So we have $x^2y \prec yz^2$ and $yz^2 \prec x^2y$,  a contradiction. So, $I$ is not a Lyubeznik ideal.
\end{ex}

\begin{rem}\label{extension} The monomial ideal $I=(x^3, x^2y, y^3, z^3)$ has only one E-minimal cover $x^2y \Box \{x^3, x^2y, y^3\}$, so it certainly is a Lyubeznik ideal by Theorem \ref{Ly ideal}. But if we add an element to the generator set to get $J=(x^3, x^2y, y^3, yz^2, z^3)$, note that $J$ is  not a Lyubeznik ideal by the example above. Hence the extension of a Lyubeznik ideal is not necessarily to be a Lyubeznik ideal. On the other hand, note that the monomial ideal $I=(x^2y^2, z^2t^2, x^2z^2, y^2t^2)$ is not a Lyubeznik ideal. But when adding $xyzt$ to its generator set, the resulting ideal $J=(x^2y^2, z^2t^2, x^2z^2, y^2t^2, xyzt)$ is a Lyubeznik ideal.
\end{rem}

\section{Some classes of Lyubeznik ideals}

As applications of Theorem 3.1 and Proposition 4.2, in this section we give several classes of Lyubeznik ideals.

 We start by observing the following two monomial ideals
 $$I=(x^2y^2, z^2t^2, x^2z^2, y^2t^2, xyzt)\supseteq J=(x^2y^2, z^2t^2, x^2z^2, y^2t^2).$$
Note that $G(I)$ contains an  element $xyzt$,
which is an inner element of every complete cover.
This is the reason why $I$ is a Lyubeznik ideal, while $J$ is not.

An element $u \in G(I)$ is called an {\it absolutely inner point} if $u$ is an inner point of every complete cover. The set of all absolutely inner points is denoted by $\mathcal{A}(G(I))$.

A monomial ideal $I$ is called a $\it cone$ ideal, if $\mathcal{A}(G(I)) \neq \emptyset$. Although the class of cone ideals are not very large, it constitutes a part of Lyubeznik ideals, as the following proposition shows.
\begin{prop}\label{cone ideal}
A cone ideal $I$ is a Lyubeznik ideal. Furthermore, the Lyubeznik resolution of $I$ is minimal under any total order $\prec$ in which the least element of $G(I)$ is an absolutely inner point.
\end{prop}

\begin{prof}
In a cone ideal $I$, let $u\in \mathcal{A}(G(I))$ be the least element of $G(I)$ in the total order $\prec$. It is easy to see that every cover $C$ of $G(I)$ is not preserved. In fact, for every cover $C$, $C\setminus \mathcal{I}(C)$ is broke by $u$, since $u\Box \overline{C}$ and $u \prec C\setminus \mathcal{I}(C)$. By Theorem \ref{Ly ideal}, $I$ is a Lyubeznik ideal under the total order $\prec$.
\end{prof}

\begin{ex}\label{ex cone 1}
The monomial ideal $I=(x^4y^4, z^4t^4, x^2y^3z^3, x^2y^2t, y^2zt^2, zyzt)$ is a cone ideal, for $xyzt$ is an inner point of every complete cover. Thus $I$ is a Lyubeznik ideal.
\end{ex}

\begin{ex}\label{ex cone 2}
Let $I$ be a monomial ideal, and let $G(I)$ be its minimal set of monomial generators. If $|G(I)|=|\mathcal{O}(G(I))|+1$, then $I$ is a cone ideal. Hence $I$ is a Lyubeznik ideal.
\end{ex}


\vs{3mm}An element $u \in G(I)$ is called a {\it c-inner point} if for every cover $C$, either $u\not\in C$ or $u$ is an inner point of $C$. We call a monomial ideal $I$ an {\it M-cone} ideal, if every M-minimal complete cover contains a c-inner point.

\begin{prop}\label{M-cone ideal}
An M-cone ideal $I$ is a Lyubeznik ideal. Furthermore, the Lyubeznik resolution of $I$ is minimal under any total order $\prec$ in which each c-inner point is  less than  the other elements of $G(I)$


\end{prop}


\begin{prof}
Under the assumption, for every cover $C$, there exists a c-inner point $u$, such that $u\Box \overline{C}$. Hence $C\setminus \mathcal{I}(C)$ is broke by $u$, since $u \prec C\setminus \mathcal{I}(C)$ by the definition of c-inner point and the construction of the total order $\prec$. Thus $C$ is not preserved. By Theorem \ref{Ly ideal}, $I$ is a Lyubeznik ideal under the total order $\prec$.
\end{prof}

\vs{3mm}The following example shows that the M-cone ideals are abundant.

\begin{ex}\label{ex M-cone}
The monomial ideal $$I=(x_1^3x_3^3x_6^3, x_2^3x_4^3x_5^3, x_1^2x_2^2, x_1x_2x_6, x_3^2x_4^2, x_3^2x_5^2, x_3x_4x_5)$$ has two M-minimal complete covers: $x_1x_2x_6 \Box \{x_1^2x_2^2, x_1^3x_3^3x_6^3, x_1x_2x_6\}$ with multidegree $x_1^3x_2^2x_3^3x_6^3$, and $x_3x_4x_5 \Box \{x_3^2x_4^2, x_3^2x_5^2, x_3x_4x_5\}$ with multidegree $x_3^2x_4^2x_5^2$. It is easy to check that each of $x_1x_2x_6$ and $x_3x_4x_5$ is a c-inner point in the corresponding M-minimal cover. Hence $I$ is an M-cone ideal, thus is a Lyubeznik ideal. Let $\prec$ be the total order defined by $$x_1x_2x_6 \prec x_3x_4x_5 \prec x_1^2x_2^2  \prec x_3^2x_4^2 \prec x_3^2x_5^2 \prec x_1^3x_3^3x_6^3 \prec x_2^3x_4^3x_5^3.$$ Then the Lyubeznik resolution of $I$ under $\prec$ is a minimal free resolution. Actually, we can choose many other total orders $\vdash$ in place of $\prec$, if only $x_1x_2x_6, x_3x_4x_5$ are less than the other elements of $G(I)$ in the order $\vdash$.
\end{ex}

\begin{rem}\label{c-inner}
Even if every M-minimal complete cover of $G(I)$ contains an inner  point,  a monomial ideal $I$ may  be a non-Lyubeznik ideal.  Example \ref{not Ly ideal} is such a counterexample.
\end{rem}

The following result on a finite partially ordered set may be a folk result. It is needed in proving Proposition \ref{tame ideal}. We include a proof for completeness. Actually,  an algorithm is provided for refining a partial order to get a total order in proving the following Lemma.

\begin{lem}\label{refine}
Every  partial order on a finite set $D$ can be refined to obtain a total order on $D$.
\end{lem}

\begin{prof}
Let $<$ be a partial order on a finite set $D$.
Note that if every element of $D$ is comparable with all the other elements, then the partial order is a total order.  Actually, we only need to show that for every element $a \in D$, there exists a partial order refined on $<$, such that $a$ is comparable with all the other elements of $D$. In fact, for a given element $a \in D$, let $A(a)$ denote  the elements less than $a$ under $<$, $B(a)$  the elements larger than $a$, $C(a)$  the elements which are not comparable with $a$. We add some relations to $<$, such that $a$ is less than every element of $C(a)$, still denoted by $<$.
 We define $u \vdash v$ if there exists a sequence $u_1, u_2, \cdots, u_k$, such that $u < u_1 < u_2 < \cdots < u_k < v$. We claim that $\vdash$ is a partial order. In fact, we only need to check that it is well defined, i.e.,  $\vdash$ fulfill the antisymmetry. If $u \vdash v$, $v \vdash u$ and $u \neq v$, then there exists sequences $u_1, u_2, \cdots, u_k$ and $v_1, v_2, \cdots, v_m$, such that $u < u_1 < u_2 < \cdots < u_k < v$ and $v < v_1 < v_2 < \cdots < v_m < u$. Hence we have $u < u_1 < u_2 < \cdots < u_k < v < v_1 < v_2 < \cdots < v_m < u$. There exists a $c \in C(a)$, such that $a < c$ is a part of the above chain, otherwise, it contradicts the assumption that the original $<$ is a partial order. So, we have a chain $u < \cdots < a < c_1 < \cdots < a < c_l < \cdots < u$ for $c_1, \cdots , c_l \in C(a)$. Thus, we have $c_l < u < a$ in the original partial order $<$,  a contradiction.
\end{prof}


\vs{3mm}If $I$ is a monomial ideal of $S=K[x_1, x_2, \cdots, x_n]$, with minimal generating set $G(I)=\{u_1, u_2, \cdots, u_s\}$ of monomials. Denote $u_i=x_1^{b_{i_1}}x_2^{b_{i_2}}\cdots x_n^{b_{i_n}}, \, i= 1, 2, \cdots, s$. Recall that $I$ is a {\it generic} monomial ideal if $b_{i_k} \neq b_{j_k}$ for all $ k=1, 2, \cdots, n$, where $i \neq j$. We call a generic monomial ideal $I$ a {\it mean ideal}, if for each pair of elements $u_i$ and $u_j$, when there exists a $k\in\{1,\cdots, n\}$, such that $0 < b_{i_k} < b_{j_k}$, then for every $k=1,\cdots, n$, either $b_{i_k} < b_{j_k}$ or $b_{j_k} = 0$. In a mean ideal, we say $u_i \prec u_j$ if there exists $k\in\{1,\cdots, n\}$, such that $0 < b_{i_k} < b_{j_k}$. A mean ideal $I$ is called a {\it tame} ideal, if $(G(I), \prec)$ is a partially ordered set.

\begin{prop}\label{tame ideal}
A tame ideal is a Lyubeznik ideal.
\end{prop}

\begin{prof}
Let $I$ be a tame ideal. By Lemma \ref{refine}, there exists a total order  $\vdash$ refined on the partial order $\prec$ on $G(I)$. By the construction of $\prec$, for every element $u\in G(I)$ and every E-minimal cover $C$ of $u$, clearly $u$ is less than every other element of $C$ in the order $\vdash$. So $C\setminus \{u\}$ is broke by $u$, and thus $C$ is not preserved. Hence $I$ is a Lyubeznik ideal by Theorem \ref{Ly ideal}.
\end{prof}

\vs{3mm}The class of tame ideals are large, as the following example shows:

\begin{ex}\label{expM-conedeal}
It is direct to check that the monomial ideals $I_1=(x^3, x^2y, y^3, y^2z, z^3)$ and $I_2=(x^4, x^3y^2, x^2yz, y^3z^2, y^4, z^3)$ are tame ideals, so they are Lyubeznik ideals by Proposition \ref{tame ideal}.
\end{ex}

\begin{rem}\label{Ly-generic} In 1998, D. Bayer and I. Peeva \cite{BPS} showed that the Scarf simplex of a generic ideal $J$ is a minimal free resolution of $J$. For a tame ideal $I$, we have proved that there exists a total order, under which the Lyubeznik resolution of $I$ is a minimal free resolution of $I$. As a kind of Lyubeznik ideal, tame ideals are generic. Note that not all the Lyubeznik ideals are generic£¬ and counterexamples appear in  Examples \ref{ex cone 1} and \ref{ex M-cone}. On the other hand, note that not all the generic ideals are Lyubeznik, e.g.,  the generic monomial ideal $(x^3, y^3, z^3, x^2y, y^2z, xz^2)$ is not a Lyubeznik ideal.
\end{rem}

\section{Examination of out points, inner points and boundary points of a complete cover}

For a given complete cover $C$, it is certainly important to distinguish the out points, inner points and boundary points in $C$.

Let $C=\{u_1, u_2, \cdots, u_m\}$, and assume $u_i=x_1^{b_{i_1}}x_2^{b_{i_2}}\cdots x_n^{b_{i_n}}, \, i= 1, 2, \cdots, m$. If $m(C)=x_1^{c_{1}}x_2^{c_{2}}\cdots x_n^{c_{n}}$, then denote

$$ d_{ij} = \left\{ \begin{array}{ccc}
0, \, \, b_{i_j} \neq c_j \\
1, \, \, b_{i_j} = c_j
\end{array}
\right., \, i=1, 2, \cdots, m; j=1, 2, \cdots, n.$$

The following proposition can be used to examine the out points, and we omit the verification here.

\begin{prop}\label{out point}
$u_i \in \mathcal{O}(C)$ if and only if
$$\prod^n_{j=1}(1-\prod_{k\neq i}(1-d_{kj}))=0.
$$\end{prop}

In order to give an algorithm to distinguish the inner points, the boundary points and the exchangeable points of $G(I)$, denote
$$B=\{j\in [n] \mid d_{ij}=0 \, for \, \forall u_i \in \mathcal{O}(C) \},\, B(u_i)=\{j\in B \mid d_{ij}=1\},$$ and
$$\mathcal{A}(u_i)=\{A \subseteq C\setminus (\mathcal{O}(C)\cup \{u_i\}) \mid \prod_{j\in B\setminus B(u_i)}(1-\prod_{k, u_k\in A}(1-d_{kj}))=1\}.$$

\begin{prop}\label{inner boundary point}
(1) $u_i \in \mathcal{I}(C)$ if and only if $$\prod_{A\in \mathcal{A}(u_i)}\prod_{j\in B(u_i)}(1-\prod_{k, u_k\in A}(1-d_{kj}))=1.$$
(2) $u_i \in \mathcal{B}(C)$ if and only if $$\prod_{A\in \mathcal{A}(u_i)}\prod_{j\in B(u_i)}(1-\prod_{k, u_k\in A}(1-d_{kj}))=0.$$
(3) For an element $u_i\in \mathcal{B}(C)$, $u_i \in \mathcal{E}(C)$ if and only if $$\prod_{A\in \mathcal{A}(u_i)}\prod_{u_{l}\in \mathcal{B}(C)\setminus A}\prod_{j\in B(u_i)}(1-\prod_{k, u_k\in A\cup \{u_{l}\}}(1-d_{kj}))=1.$$
\end{prop}

The  above propositions follow from  explicit computations, so we omit  detailed verifications too.

\end{document}